\def\Date{July 27, 2009}
\magnification=\magstephalf
\frenchspacing\parskip4pt plus 1pt
\advance\vsize6mm

\newread\aux\immediate\openin\aux=\jobname.aux 
\ifeof\aux\message{ <<< Run TeX a second time >>> }
\batchmode\else\input\jobname.aux\fi\closein\aux
\input \jobname.def\immediate\openout\aux=\jobname.aux

\draftfalse
\vglue1cm

{\parindent0pt\footnote{}{\ninerm Mathematics Subject
Classification (2000): 32V20, 17B66, 16W10, 32M15}}

\headline={\ifnum\pageno>1\eightrm\ifodd\pageno\hfill~~ CR-quadrics
with a symmetry property \hfill{\tenbf\folio}\else{\tenbf\folio}\hfill
W. Kaup\hfill\fi\else\hss\fi}

\centerline{\Gross CR-quadrics with a symmetry property}
\medskip\bigskip
\centerline{{WILHELM KAUP}}

\KAP{Introduction}{Introduction}
A well studied class of CR-submanifolds in a complex linear space
consists of
the quadrics, that is, of real quadratic submanifolds of the form
$$Q=\{(w,z)\in\CC^{k+n}:w+\overline w=h(z,z)\}\,,$$ where
$h:\CC^{n}\times\CC^{n}\to\CC^{k}$ is a non-degenerate hermitian form
such that the image of $h$ spans all of $\CC^{k}$. It is well known
and easy to see that the group $\Aff(Q)$ of all (complex affine)
transformations leaving $Q$ invariant acts transitively on
$Q$. Besides these global CR-automorphisms of $Q$, in general there
exist also non-affine {\sl local} CR-automorphisms (between domains)
of $Q$, which cannot be extended to global CR-automorphisms of $Q$. By
$\Lit{BAJT}$ every such local (smooth) CR-isomorphism is real-analytic
and by \Lit{TUMA} extends to a birational transformation of
$\CC^{k+n}$. All birational transformations obtained this way generate
a group that we denote by $\Bir(Q)$ in the following. It is by no
means evident, but shown in \Lit{ISKA}, that
$g\big(Q\cap\reg(g)\big)=Q\cap\reg(g^{-1})$ holds, where
$\reg(g)\subset\CC^{k+n}$ is the Zariski open subset of all regular
points of the birational transformation $g$ with $Q\cap\reg(g)$ being
a dense domain in $Q$.

Up to CR-isomorphism there exists a unique homogeneous real-analytic
CR-manifold $\hat Q$ containing $Q$ as open dense CR-submanifold in
such a way that every CR-isomorphism between domains in $Q$ extends to
a global CR-automorphism of $\hat Q$. In particular, the group
$\Bir(Q)$ can be canonically identified with the CR-automorphism group
$\Aut(\hat Q)$ of the {\sl extended quadric} $\hat Q$. The group
$G:=\Aut(\hat Q)$ has no center and can be realized as closed subgroup
$G\subset\SL_{N}(\CC)$ for some integer $N\ge2$ is such a way that
$\hat Q$ is a $G$-orbit in the complex projective space
$Z:=\PP(\CC^{N})$. In fact, there are complex-algebraic subvarieties
$A\subset B\subset Z$ such that the $G$-orbit $\hat Q$ is a closed
CR-submanifold of $Z\backslash A$ and $Q=\hat Q\backslash B$ (for all
this and further details compare \Lit{ISKA}). In particular, the group
$G$ inherits a Lie group structure from $\SL_{N}(\CC)$ and with it is
a transitive real-analytic transformation group on $\hat Q$. It can be
shown that $G$ has only finitely many connected components.

\medskip A convenient method for the study of a Lie group action is
given by the associated infinitesimal action in terms of vector
fields.  For every domain $D\subset Q$ let $\hol(D)$ be the real Lie
algebra of all real-analytic infinitesimal CR-transformations on $D$
(that is of all real-analytic vector fields on $D$ whose local flows
consist of CR-transformations). By \Lit{BELO} every vector field in
$\hol(D)$ extends to a (complex) polynomial vector field of degree
$\le2$ on $\CC^{k+n}$.  This implies, in particular, that
$\7g:=\hol(Q)$ has finite dimension and that for every domain
$D\subset Q$ the restriction operator $\7g\to\hol(D)$ is an
isomorphism of real Lie algebras. Furthermore, $\7g$ is in a canonical
way the Lie algebra of the Lie group $G=\Aut(\hat Q)$. The vector
field $\zeta:=2w\dd w+z\dd z\in\7g$ gives a canonical grading
$$\7g=\7g^{-2}\oplus\7g^{-1}\oplus\7g^{0}\oplus\7g^{1}\oplus
\7g^{2}\Leqno{AL}$$ into the $\ad(\zeta)$-eigenspaces, where the Lie
algebra $\aff(Q):=\7g^{-2}\oplus\7g^{-1}\oplus\7g^{0}\subset\7g$ is
the Lie subalgebra of the subgroup $\Aff(Q)\subset G$.  The affine
subalgebra $\aff(Q)$ has an explicit description in terms of the
hermitian form $h$, see \Ruf{DN}, whereas it seems to be unknown how
big the nilpotent Lie subalgebra $\7g^{+}:=\7g^{1}\oplus\7g^{2}$ can
be in terms of $k,n$ (it can definitely be zero). Once $\7g^{+}$ is
explicitly known also the nilpotent closed subgroup
$G^{+}:=\exp(\7g^{+})\subset G$ is explicitly known - indeed,
$\7g^{+}$ is $\ad$-nilpotent so that $\exp:\7g^{+}\to G^{+}$ is a
polynomial homeomorphism. The group $G$ is generated by the connected
subset $\exp(\7g)$ together with the linear subgroup
$\GL(Q)\subset\Aff(Q)$. The Lie algebras $\7g^{k}$ and hence $\7g$
itself can be explicitly determined by solving certain linear
equations, a comparatively much easier task than the explicit
determination of the linear group $\GL(Q)$, where high order
polynomial equations have to be solved.

\medskip The intention of this short note is to discuss several
classes of examples for $Q$ as above where
$\7g^{+}=\7g^{1}\oplus\7g^{2}$ is `big' and has a simple description.
For that we introduce the Property (S): There exists a transformation
$\gamma=\gamma^{-1}\in G=\Aut(\hat Q)\cong\Bir(Q)$ such that
$\Ad(\gamma)(\zeta)=-\zeta$ for the adjoint representation
$\Ad:G\to\Aut(\7g)$. Then, if such a $\gamma$ is obtained in a
concrete situation (by guess-work or any other form of computation)
the obvious formula $\7g^{k}=\Ad(\gamma)(\7g^{-k})$ allows to
immediately write down $\7g^{+}=\7g^{1}\oplus\7g^{2}$ explicitly.  The
same holds with $G^{+}=\gamma H\gamma$, where
$H:=\exp(\7g^{-2}\oplus\7g^{-1})$ is the Heisenberg group. As
indicated above, the determination of $\GL(Q)$ is more involved.

\bigskip The paper is organized as follows: After the necessary
preliminaries in Section \ruf{Preliminaries} we introduce in Section
\ruf{symmetry} the symmetry Property (S) and present with Example
\ruf{EX} our basic class of quadrics having this property.  In Section
\ruf{Tensoring} we obtain for every $Q$ from the basic class by
tensoring with an arbitrary unital (associative) complex \sta{} $A$ of
finite dimension a new quadric $Q(A)$ that also has Property (S).
In the final section we briefly explain how from the classification of
irreducible bounded symmetric domains of non-tube type further
quadrics with Property (S) can be obtained.

\KAP{Preliminaries}{Preliminaries}

Let $W,Z$ be complex vector spaces of finite positive
dimension. Suppose that on $W$ a conjugation $w\mapsto\overline w$ is
given and put $V:=\{w\in W:\overline w=w\}$. Then for every
sesquilinear form $h:Z\times Z\to W$ (complex linear in the first and
conjugate linear in the second variable) the real-algebraic subset
$$Q=Q_{h}=\{(w,z)\in W\times Z:w+\overline w=h(z,z)\}\Leqno{AS}$$ is
called a {\sl standard} quadric in the following, provided \0
$h(z,z')=0$ for all $z'\in Z$ implies $z=0\,$ (non-degeneracy), \1 $V$
is the linear span over $\RR$ of all vectors $h(z,z)$, $z\in Z$
(minimality).

\bigskip It is clear that $Q$ is invariant under the two 1-parameter
groups of linear transformations
$$(w,z)\mapsto(e^{2t}w,e^{t}z)\steil{and}(w,z)\mapsto(w,e^{it}z)
\,,\quad t\in \RR.$$ Therefore $\7g:=\hol(Q)$ contains the commuting
linear vector fields
$$\zeta:=2w\dd w+z\dd z\Steil{and}\chi:=\,iz\dd z\,,\Leqno{DL}$$ and
$\7g\!+i\7g$ contains the Euler field
$\eta=(\zeta-i\chi)/2$. 

\medskip Denote by $\7P$ the complex Lie algebra of all (complex)
polynomial vector fields on $E:=W\oplus Z$. The vector field
$\zeta\in\7g\subset\7P$ induces a $\2Z$-grading
$$\7P=\bigoplus_{k\in\2Z}\7P^{k}\,,\quad[\7P^{k},\7P^{\ell}]\subset
\7P^{k+\ell}\steil{for}\7P^{k}:=\{\xi\in\7P:[\zeta,\xi]=
k\xi\}\Leqno{BF}$$ of $\7P$ with $\7P^{k}=0$ if $k<-2$ and induces
also the grading \Ruf{AL} with $\7g^{k}:=\7g\cap\7P^{k}$, compare
\Lit{BELO}. The subalgebra $\7g^{-2}\oplus\7g^{0}\oplus\7g^{2}$ is the
kernel of $\ad(\chi)$, while the restriction of $\ad(\chi)$ to the
invariant subspaces $\7g^{-1}$, $\7g^{1}$ has the eigenvalues $\pm i$.
The following is well known and easily verified:
$$\eqalign{\7g^{-2}&=\big\{a\dd w:a\in
iV\big\},\cr\7g^{-1}&=\big\{h(z,c)\dd w+c\dd z:c\in
Z\big\},\cr\7g^{0}\mskip11mu&=\big\{aw\dd w+bz\dd
z:a\in\gl(V),b\in\gl(Z)\steil{with}ah(z,z)=
h(bz,z)+h(z,bz)\big\}\;.}\Leqno{DN}$$ \noindent The derived algebra
$\7d:=[\7g,\7g]$ is an ideal in $\7g$ and has the grading
$\7d=\bigoplus_{|k|\le2}\7d^{k}$, where
$$\7d^{k}=\cases{\7g^{k}&$k\ne0$\cr[\7g^{-2},\7g^{ 2}]+ [\7g^{-1},\7g^{
1}]+[\7g^{ 0},\7g^{ 0}]&$k=0\;$.\cr}$$ Furthermore,
$\7g^{-}:=\7g^{-2}\oplus\7g^{-1}$ and $\7g^{+}:=\7g^{1}\oplus\7g^{2}$
are nilpotent Lie algebras of step 2.

\medskip Let $\hat Q$ be the extended quadric and $G:=\Aut(\hat Q)$,
compare Section \ruf{Introduction}. Then $\7g=\hol(Q)$ can also be
identified canonically with $\hol(\hat Q)$. Consider the following
subgroups of $G$
$$\eqalign{G_{0}:&=\{g\in G:g(0)=0\}\cr G^{\pm}:&=\exp(\7g^{\pm})\cr
\GL(Q):&=\{g\in\GL(E):g(Q)=Q\}\cr
\Aff(Q):&=\{g\in\Aff(E):g(Q)=Q\}=\GL(Q)\,\ltimes\, G^{-}\cr }$$ with
$\Aff(E)$ the group of all complex affine automorphisms of
$E$. Then$$\GL(Q)=\{(f\times g)\in\GL(V)\times\GL(Z)\subset\GL(E):
fh(z,z)=h(gz,gz)\}\Leqno {UV}
$$ is a real algebraic subgroup of $\GL(E)$ with Lie algebra
$\gl(Q):=\7g^{0}\subset\gl(E)$ and, in particular, has only finitely
many connected components. Every $G^{\pm}$ is a connected nilpotent
closed subgroup of $G$ with Lie algebra $\7g^{\pm}$. For instance,
$G^{-}=\exp(\7g^{-2})\exp(\7g^{-1})$ is the group of all affine
transformations of the form
$$(w,z)\mapsto(w+a+h(z,b),z+b)\,,\quad (a,b)\in Q\,,\Leqno{DK}$$ which
acts simply transitively on $Q$ and is called the {\sl Heisenberg
group}.

\medskip Every $g\in G=\Aut(\hat Q)$ acts on its Lie algebra $\7g$ by
$\Ad(g)\in\Aut(\7g)$, here given in terms of vector fields by
$$\Ad(g)\big(f(z)\dd z\big)=h(z)\dd
z\Steil{with}h(g(z))=g'(z)(f(z))\,,$$ where $g'(z)\in\End(E)$ is the
derivative of $g$ at $z$. The group $G$ has no center, that is, the
group homomorphism $\Ad:G\to\Aut(\7g)$ is injective.

\medskip The following result will not be used later but may be of
independent interest.

\Lemma{} For every standard quadric $Q$ the extended quadric $\hat Q$
is simply connected.

\Proof Denote by $\pi:\tilde Q\to\hat Q$ the universal covering of
$\hat Q$.  Then by \Lit{ISKA} there exists a complex manifold $X$,
containing $\hat Q$ as generic real-analytic submanifold, together
with a complex-analytic subset $A\subset X$ such that $Q=\hat
Q\backslash A$. There exists a complex manifold $\tilde X$,
containing $\tilde Q$ as generic real-analytic submanifold, in such a
way that $\pi$ extends to a holomorphic map $\pi:\tilde X\to X$. This
implies that $\tilde A:=\pi^{-1}(A)$ is complex-analytic in $\tilde
X$ and hence that $\pi^{-1}(Q)=\tilde Q\backslash\tilde A$ is
connected by Lemma 2.2 in \Lit{FKAP}.  Since $Q$ is simply connected the
covering map $\pi:\tilde Q\to\hat Q$ must be a homeomorphism.\qed

\KAP{symmetry}{The symmetry property}

With the notation of Section \ruf{Preliminaries} fix a standard
quadric $Q$ and consider the following symmetry property:

\kap{Property (S)} There exists an automorphism $\gamma=\gamma^{-1}\in
G=\Aut(\hat Q)$ with $\Ad(\gamma)(\zeta)=-\zeta$.\nline We call
$\gamma$ also a {\sl symmetry of the quadric} $Q$. There may not exist
a fixed point of $\gamma$ in the extended quadric $\hat Q$, insofar
$\gamma$ is not necessarily a CR-{\sl symmetry} of $\hat Q$ in  the
sense of \Lit{KAZA}.

\bigskip If the symmetry property (S) is satisfied with $\gamma$ then
$\Ad(\gamma)\in\Aut(\7g)$ permutes the eigenspaces of $\ad(\zeta)$ in
$\7g$, more precisely,
$$\Ad(\gamma)(\7g^{k})=\7g^{-k}\steil{for all}k,\steil{in
particular,}\Leqno{DM}$$
$$\dim(\7g^{k})=\dim(\7g^{-k})\steil{for
all}k\Steil{and}[\7g^{1},\7g^{1}]=\7g^{2}\,.$$ The symmetry $\gamma$
is not uniquely determined, every $g\gamma g^{-1}$ with $g\in\GL(Q)$
is also a symmetry of $Q$. Using \Ruf{DM} the spaces $\7g^{1}$ and
$\7g^{2}$ can be explicitly computed from \Ruf{DN}. In Section
\ruf{Tensoring} we given an example of this method. This also works on
the group level. Indeed, the inner automorphism $\Int(\gamma)$ of $G$
defined by $g\mapsto \gamma g\gamma$ (note that $\gamma^{-1}=\gamma$
by definition) satisfies $$G^{+}=\gamma\,
G^{-}\gamma\Steil{and}G_{0}=\GL(Q)\,\ltimes\,G^{+}=
\gamma\,\Aff(Q)\,\gamma\,.\Leqno{DO}$$

As a consequence we state
\Proposition{} The group $G$ is generated by the subgroup $\Aff(Q)$
and $\gamma$.\Formend

\bigskip In the following we give some examples of standard quadrics
having a symmetry.  We start with the case of hyperquadrics.

\Example{EX} Suppose that in $\CC^{n+1}$ with coordinates
$(w,z_{1},\dots,z_{n})$ the quadric $Q$ is given by
$$Q=\Big\{(w,z)\in\CC^{n+1}:w+\overline w=\sum_{1\le j\le
k}|z_{j}|^{2}-\sum_{k<j\le n}|z_{j}|^{2}\Big\}\;,$$ where $0\le k\le
n$ is a fixed integer. Then $(w,z)\mapsto (w^{-1},w^{-1}z)$ defines a
symmetry $\gamma$ of $Q$. Obviously, there is a fixed point $(1,z)\in
Q$ of $\gamma$, provided $k>0$ (in case $k=0$ the symmetry
$(w,z)\mapsto (w^{-1},-w^{-1}z)$ would have a fixed point $(-1,z)\in
Q$). The Lie algebra $\hol(Q)$ is isomorphic to $\su(p,q)$ with
$p=k{+}1$ and $q=n{+}1{-}k$. Clearly, replacing $k$ by $n-k$ gives a
linearly equivalent quadric. \Formend

\medskip This example can be generalized to higher codimensions. For
every matrix $w$ we denote by
$w^{*}$ its adjoint (conjugate transpose).

\Example{EY} Let $m,n\ge1$ be fixed integers. Then, for every
hermitian matrix $\beta\in\GL_{n}(\CC)$,
$$Q:=\big\{(w,z)\in\CC^{m\times m}\times\CC^{m\times n}:w+w^{*}=z\beta
z^{*}\big\}$$ is a standard quadric with CR-codimension $m^{2}$, and
$(w,z)\mapsto (w^{-1},w^{-1}z)$ defines a symmetry $\gamma$ of $Q$. In
case $n<m$ there is no fixed point of $\gamma$ in $\hat Q$. On the
other hand, in case $n\ge m$ there exists a fixed point $(\One,z)\in
Q$, provided $k\ge m$ for the number $k$ of positive eigenvalues of
$\beta.$ The Lie algebra $\hol(Q)$ is isomorphic to $\su(p,q)$ with
$p=k+m$ and $q=n+m-k$. The extended quadric $\hat Q$ can be realized
as follows: Choose on $\CC^{2m+n}$ a non-degenerate hermitian form
$\psi$ of type $(p,q)$ and denote by $\2G$ the Grassmannian of all
linear $m$-spaces in $\CC^{p+q}$. Then the compact real-analytic
submanifold $S:=\{L\in\2G:\psi_{|L}=0\}$ of $\2G$ is CR-isomorphic to
$\hat Q$. In fact, $S$ is the unique closed orbit of the unitary group
$\SU(\psi)\cong\SU(p,q)$ acting on the Grassmannian. For $m=1$ we get
back Example \ruf{EX} (up to an affine transformation) and $\2G$ is
the complex projective space $\PP(\CC^{n{+2}})$.

\medskip It is obvious that for every pair $Q'$, $Q''$ of standard
quadrics with Property (S) also the direct product $Q:=Q'\times Q''$
is a standard quadric with Property (S).  In the next section we
describe a more interesting method to produce new standard quadrics
out of known ones.

\KAP{Tensoring}{Tensoring with \sta s}

In the following we use the notion of a \sta, that is, a
complex associative algebra $A$ with product $(w,z)\mapsto wz$ and
(conjugate linear) involution $z\mapsto z^\star$ satisfying
$z^{\star\star}=z$ and $(zw)^\star=w^\star\! z^\star$ for all $w,z\in
A$. Then the {\sl selfadjoint part} $A_{sa}:=\{z\in A:z^\star=z\}$ is
a real {\sl Jordan subalgebra} with respect to the anti-commutator
product $x\circ y:=(xy+yx)/2$ and $iA_{sa}$ is a real {\sl Lie
subalgebra} with respect to the commutator product
$[x,y]=xy-yx$. Clearly $A=A_{sa}\oplus\,iA_{sa}$, and $A_{sa}$ is an
associative subalgebra of $A$ only if $A$ is commutative. 

Here we assume without further notice that every \sta{} has finite
dimension.  In addition we also assume that $A$ has a unit $e$, which
then is contained in $A_{sa}$. The subgroup $\G(A)\subset A$ of all
invertible elements is Zariski open in $A$ and a connected complex Lie
group. In particular, $\G(A)\subset\GL(A)$ is generated by the image
of the exponential map $\exp:A\to\G(A)$. In case $A$ is commutative,
$\exp$ is a surjective group homomorphism as well as a locally
biholomorphic covering map.

Denote by
$$\Aut(A,\star):=\big\{g\in\GL(A):g(ac^{\star})=g(a)g(c)^{\star}\big\}$$
the \sta{} automorphism group of $A$. Then $\Aut(A,\star)$ is a real
linear algebraic subgroup of $\GL(A)$ leaving $A_{sa}$ invariant and
has Lie algebra
$$\der(A,\star):=\{\delta\in\gl(A):\delta(ac^{\star})=
\delta(a)c^{\star}+a\delta(c)^{\star}\}\,.$$ In case that $A$ is commutative,
$\Aut(A,\star)$ can be identified with the real algebra automorphism
group of the associative algebra $A_{sa}$ and $\der(A,\star)$ can be
identified with the derivation algebra of $A_{sa}$.

For every integer $m\ge 1$ the matrix algebra $M:=\CC^{m\times m}$ is
a unital \sta{} with respect to the usual adjoint $*$ as involution,
and $M$ is even simple as complex algebra. Every other involution
$\star$ on $M$ is of the form $z^{\star}=\alpha z^{*}\alpha^{-1}$ for some
$\alpha=\alpha^{*}\in\GL_{m}(\CC)$. Also, the product algebra $M\times M$
becomes a simple \sta{} with respect to the involution
$(a,b)^{\star}:=(b^{*},a^{*})$.  By Wedderburn's Theorem \Lit{LANG}
every semi-simple complex unital algebra $A$ (of finite dimension) is
a unique direct product $A=\prod_{j\in J} A_{j}$ of simple unital
algebras $A_{j}\cong\CC^{m_{j}\times m_{j}}$. To every involution
$\star$ of $A$, making it to a \sta, there is an involution $j\mapsto
j^{\bullet}$ of the index set $J$ such that
$(A_{j})^{\star}=A_{j^{\bull}}$ for every $j\in J$. Then, choosing a
minimal subset $K\subset J$ with $J=K\cup K^{\bullet}$, we get the
representation $A=\prod_{j\in K} (A_{j}+A_{j^{\bull}})$ as direct
product of simple \sta s.

In general the algebra $A$ has a radical (the intersection of all
maximal left ideals in $A$). Every involution of $A$ leaves $\Rad(A)$
invariant and makes the semi-simple algebra $A/\Rad(A)$ to a \sta.  It
can be shown that even in the commutative case there are uncountably
many isomorphy classes of (finite dimensional) \sta s.

\bigskip Now consider for given $W,Z,h$ the standard quadric $Q=Q_{h}$
defined as in \Ruf{AS} and fix a unital \sta{} $A$. Put
$V:=\{w\in W:\overline w=w\}$ and also
$$\3W:=W\otimes_{\CC}\!A,\quad\3Z:=W\otimes_{\CC}\!A\steil{and}\3V:
=V\otimes_{\RR}\!A_{sa}\,.$$ Then there exists a unique conjugation
$\3w\mapsto\overline\3w$ on $\3W$ such that $\overline\3w=\overline
w\otimes a^{*}$ for every $\3w=w\otimes a\in\3W$. The fixed point set
of this conjugation is $\3V$, considered in a canonical way as
$\RR$-linear subspace of $\3W$. Also there exists a unique hermitian
form $$\3h:\3Z\times\3Z\to\3W\steil{satisfying}\3h(x\otimes a,y\otimes
b)=h(x,y)\otimes ab^{\star}$$ defining the quadric
$$\3Q:=Q(A):=\{(\3w,\3z)\in \3W\times \3Z:\3w+\overline
\3w=\3h(\3z,\3z)\}\,.$$ If we denote by $e\in A$ its unit we realize
$V,W,Z$ as linear subspaces of $\3V,\3W,\3Z$ by identifying every
$w\in W$, $z\in Z$ with $w\otimes e$, $z\otimes e$ respectively. In
this sense $Q=\3Q\cap(W\times Z)\subset\3W\times\3Z$. Furthermore, for
every pair $A$, $B$ of unital \sta s we have $Q(A\times B)\cong
Q(A)\times Q(B)$ and $Q(A)(B)\cong Q(A\otimes B)$, where product and
involution on $A\otimes B$ are uniquely determined by $(a\otimes
b)(c\otimes d)=(ac\otimes bd)$ and $(a\otimes
b)^{\star}=(a^{\star}\otimes b^{\star})$. 
\medskip
\kap{Tensored Example {\EX}:} We consider Example {\EX} in the lowest
possible dimension $n=1$, that is without loss of generality,
$$Q=\{(w,z)\in\CC^{2}:w+\overline w=z\overline z\}$$ is the Heisenberg
sphere in dimension $2$. With $A$ a fixed \sta{} of complex dimension
$d$ then
$$\3Q=Q(A)=\{(w,z)\in A^{2}:w+w^{\star}=zz^{\star}\}\Leqno{BR}$$ is a
standard quadric of CR-codimension $d$, and a symmetry is given by
$\gamma(w,z)=(w^{-1},w^{-1}z)$.

\Proposition{DQ} For $\3Q$ in \Ruf{BR} and
$\hol(\3Q)=\7g^{-2}\oplus\7g^{-1}\oplus\7g^{0}\oplus\7g^{1}\oplus
\7g^{2}$ we have
$$\eqalign{\7g^{-2}&=\big\{a\dd{w}:a\in iA_{sa}\big\}\cr
\7g^{-1}&=\big\{zc^{\star}\dd{w}+c\dd{z}:c\in A\big\}\cr
\7g^{0}&=\big\{(aw+wa^{\star})\dd{w}+az\dd{z}:a\in
A\big\}\;\rtimes\;\big\{\delta(w)\dd{w}+
\delta(z)\dd{z}:\delta\in\der(A,\star)\big\}\cr
\7g^{1}&=\big\{zc^{\star}w\dd{w}+(zc^{\star}z+wc)\dd{z}:c\in
A\big\}\cr \7g^{2}&=\big\{waw\dd{w}+waz\dd{z}:a\in
iA_{sa}\big\} \;.\cr}$$

\Proof $\7g^{k}$ for $k<0$ is \Ruf{DN} and for $k>0$ follows
immediately by applying \Ruf{DM} to $\7g^{-1}$, $\7g^{-2}$. For $k=0$
the claim is a direct consequence of the following proposition.\qed

\Proposition{DP} For $\3Q$ in \Ruf{BR} and the subgroups $\GL(\3Q)$,
$G^{+}$ of $G=\Aut(\hat{\3Q})$ as defined in Section
\ruf{Preliminaries} we have, where $e\in A$ is the unit:
$$\eqalign{\GL(\3Q)
&=\big\{(w,z)\mapsto(awa^{\star},az):a\in\G(A)\big\}\;\rtimes\;
\big\{g\times g:g\in\Aut(A,\star)\big\}\cr
G^{+}&=\{(w,z)\mapsto(e+wa+zc^{\star})^{-1}
(w,z+wc):(a,c)\in \3Q\} \;.\cr}$$ 

\Proof The second equation follows immediately with \Ruf{DO} applied
to \Ruf{DK}.\nline Next consider an arbitrary $\phi\in\GL(\3Q)$. Then
$\phi=f\times g\;\in\; \GL(A_{sa})\times\GL(A)$ with
$f(xy^{\star})=(gx)(gy)^{\star}$ for all $x,y\in A$ by $\Ruf{UV}$. We
claim that $g(e)$ is invertible in $A$ (compare also Lemma 2.5 in
\Lit{ESSC}): Indeed, for $x=e$ we have $f(y)=(ge)(gy)^{\star}$, and
choosing $y$ with $f(y)=e$ we find $e=(ge)(gy)^{\star}$, that is,
$g(e)\in \G(A)$.\nline Obviously the group
$K:=\{(w,z)\mapsto(awa^{\star},az):a\in \G(A)\}$ is contained in
$\GL(Q)$ and the group $\{\beta\in GL(A):(\alpha\times \beta)\in
K\steil{for some}\alpha\in\GL(A_{sa})\}$ acts simply transitively on
$\G(A)$.  We may therefore assume that $g(e)=e$ for $\phi=f\times g$.
For every $x\in A_{sa}$ then
$f(x)=f(xe^{\star})=(gx)(ge)^{\star}=g(x)$, that is $f=g$. Finally,
$f(xy^{\star})=(gx)(gy)^{\star}$ implies $g\in\Aut(A,\star)$.  This
implies that $K$ is normal in $\GL(\3Q)$ and that $\GL(\3Q)$ is a
semi-direct product of the claimed form. \qed

\bigskip

In \Lit{ESSC} the notion of a {\sl Real Associative Quadric} (RAQ for
short) has been introduced and for every quadric $Q$ of this type the
subgroups $\GL(Q)$ and $G^{+}$ of $\Aut(Q)$ have been explicitly
described. In our terminology, the RAQs are just the tensored quadrics
$Q(A)$ of type \Ruf{BR}, where $Q\subset\CC^{2}$ is the Heisenberg
quadric in dimension 2 and $A$ is a {\sl commutative} \sta.  In this
special case of a commutative \sta{} $A$ the group $\Aut(A,\star)$
obviously can also be written in the following way:
$$\Aut(A,\star)=\Aut(R)\Steil{ for the commutative real
algebra}R=A_{sa}\,.$$ For every $a\in R$ we denote by $L(a)\in\
\End(R)$ the left multiplication $x\mapsto ax$ and identify the
algebra $R$ with its image $L(R)\subset\End(R)$. Then the normalizer
$\{g\in\GL(R):g\,L(R)=L(R)\,g\}$ is canonically isomorphic to
$\Aut(R)$.  Via this identification, for the commutative case the
descriptions of $\GL(\3Q)$ and $G^{+}$ in Proposition \ruf{DP} occur
already in \Lit{ESSC}.

\medskip \kap{Tensored Example {\EY}:} For $m,n\ge1$ and hermitian
matrix $\beta\in\GL_{n}(\CC)$ as in Example
{\EY} we consider the tensored quadric $\3Q=Q(A)$, where $A$ is an
arbitrary \sta, that is,
$$\3Q=\big\{(w,z)\in A^{m\times m}\times A^{m\times
n}:w+w^{\star}=z\beta z^{*}\big\}\,,$$ where for every matrix
$a=(a_{jk})\in A^{r\times s}$ the adjoint $a^{*}\in A^{s\times r}$ is
the matrix $(a_{kj}^{\star})$. Then $\3Q$ is a standard quadric and,
as before, $(w,z)\mapsto(w^{-1},w^{-1}z)$ defines a symmetry. It is
seen immediately that for $G^{+}$ and $\7g=\hol(\3Q)$ we have
$$\eqalign{G^{+}&= \{(w,z)\mapsto(e+wa+z\beta c^{\star})^{-1}
(w,z+wc):(a,c)\in \3Q\} \;.\cr \7g^{-2}&=\big\{a\dd{w}:a\in A^{m\times
m},\,a+a^{\star}=0\}\cr \7g^{-1}&=\big\{z\beta
c^{\star}\dd{w}+c\dd{z}:c\in A^{m\times
n}\big\}\cr\7g^{1}&=\big\{z\beta c^{\star}w\dd{w}+(z\beta
c^{\star}z+wc)\dd{z}:c\in A^{m\times n}\big\} \cr \7g^{2}&=
\big\{waw\dd{w}+waz\dd{z}:a\in A^{m\times
m},\,a+a^{\star}=0\big\}\;,\cr }$$ and a simple check reveals
$[\7g^{j},\7g^{k}]=\7g^{j+k}$ for all $j,k\in\{0,\pm1,\pm2\}$ with
$j+k\ne0$.

More involved is the linear group $\GL(\3Q)$.  Notice that the group
$\GL_{n}(A):=\G(A^{n\times n})$ acts by matrix multiplication from the
right on $A^{m\times n}$. A subgroup is the $\beta$-{\sl unitary
group}
$$\U_{\beta}(A):=\{u\in\GL_{n}(A):u\beta u^{\star}=\beta\}\steil{with
Lie algebra}\7u_{\beta}(A):=\{u\in A^{n\times n}:u\beta+\beta
u^{\star}=0\}\,.$$ 
For every matrix space
$A^{r\times s}$ we have the embedding $\End(A) \into\End(A^{r\times
s})$ given by $f\cd a=(f(a_{jk}))$ for every $a=(a_{jk})\in A^{r\times
s}$.
Then a simple computation gives that
$$\B:=\big\{(w,z)\mapsto(a(g\cd w)a^{\star},a(g\cd
z)u):a\in\GL_{m}(A),\; u\in\U_{\beta}(A),\;g\in\Aut(A,\star)\big\}$$
is a subgroup of $\GL(\3Q)$.  In case $A$ is a semi-simple, up to a
permutation of its simple factors, every element of $\Aut(A,\star)$ is
an inner $\star$-automorphism of the form $a\mapsto uau^{\star}$ with
$u^{\star}=u^{-1}\in \G(A)$. In particular, $\Aut(\CC,\star)$ is the
trivial group and it can be seen that $\B=\GL(Q)$ holds in case
$A=\CC$.  Also $\B=\GL(\3Q)$ in case $m=n=1$ by
Proposition \ruf{DP}.

\bigskip We discuss briefly another local realization of the above
tensored quadrics $\3Q=Q(A)$: Fix integers $r>m\ge1$ and a hermitian
matrix $\alpha\in\GL_{r}(\CC)$ having at least $m$ positive
eigenvalues. For fixed \sta{} $A$ put $E:=A^{m\times r}$ and let
$\One\in\GL_{m}(A)$ be the unit matrix. Then we call
$$\3S:=\{z\in E:z\alpha z^{\star}=\One\}$$ a {\sl generalized sphere}.
Without loss of generality we may assume that
$\alpha=\One\times\beta\in\GL_{m}(\CC)\times\GL_{n}(\CC)$ for
$n:=r-m$ and some hermitian matrix $\beta\in\GL_{n}(\CC)$. Now
put $\;W:=A^{m\times m}$, $\;Z:=A^{m\times n}$.  Then
$E=W\oplus Z$ in a canonical way and the quadric $\3Q:=\{(x,y)\in
W\times Z:x+x^{*}=y\beta y^{*}\}$ is locally CR-isomorphic to $\3S$.
Indeed, consider the Cayley transformation $\kappa\in\Bir(E)$ defined
on $E$ by
$$\kappa(x,y)=(\One-x)^{-1}(\One+x,\sqrt2y)\,.\Leqno{UQ}$$ Then $\kappa^{-1}(
x,y)=(x+\One)^{-1}(x-\One,\sqrt2y)\,$ and a simple computation shows
that the birational transformation $\kappa$ gives a CR-isomorphism
$$\kappa:\3S\cap\reg(\kappa)\;\to\; \3Q\cap\reg(\kappa^{-1})\,.$$
Since $\3Q\cap\reg(\kappa^{-1})$ is connected by Lemma 2.2 in
\Lit{FKAP}, also $S\cap\reg(\kappa)$ is connected. Consequently $\3S$
is a connected generic real-analytic CR-submanifold of $E$.
Furthermore, $\kappa$ induces an isomorphism between the real Lie
algebras $\7g=\hol(\3Q)$ and $\7s:=\hol(\3S)$. Since
$\kappa=\exp(\xi)$ for some $\xi\in\7l:=\7g+i\7g$, we have also
$\7l=\7s+i\7s$ and every vector field in $\7s$ is polynomial of degree
$\le2$ on $E$. For every $a\in E$ application of the vector field
$\xi:=(a-z\alpha a^{\star}z)\dd z$ on $E$ to the defining equation for
$\3S$ gives $$\eqalign{\xi(z\alpha z^{\star}-\One)&=(a-z\alpha
a^{\star}z)\alpha z^{\star}+z\alpha(a^{\star}-z^{\star}a\alpha
z^{\star})\cr &=(\One-z\alpha z^{\star})a\alpha z^{\star}\;+\;z\alpha
a^{\star}(\One-z\alpha z^{\star})\,.\cr}$$ This shows that $\xi$ is
tangent to $\3S$ and we get the decomposition
$$\7s=\7k\oplus\7p\Steil{with}\7k:=\{\xi\in\7s:\xi_{0}=0\}\steil{and}
\7p:=\{(a-z\alpha a^{\star}z)\dd z:a\in E\} \,.$$ Notice that the
evaluation map $\epsilon_{0}:\7s\to E$ at the origin induces an
$\RR$-linear isomorphism of $\7p$ onto $E$. The Lie subalgebra $\7k$
contains the multiple $\delta:=iz\dd z$ of the Euler field, and $\7k$
resp. $\7p$ are the $0$- resp. $-1$-eigenspaces of $(\ad\,\delta)^{2}$
in $\7s$.  Also, every vector field in $\7k$ is linear, and
$\7k=\gl(\3S)$ is the Lie algebra of the linear algebraic group
$\K:=\GL(\3S)=\{g\in\GL(E):g(\3S)=\3S\}$. It is clear that the group
$\U:=\U_{\one}\times\U_{\alpha}\subset\GL_{m}(A)\times\GL_{r}(A)$ acts
linearly on $\3S$ via $z\mapsto uzv^{\star}$, and it can be seen that
there is an open orbit in $\3S$ for this action .

\KAP{Connection}{Links with bounded symmetric domains}

For the special case $\beta=\One$ in Example \ruf{EY} the
quadric $Q$ coincides with the \v Silov boundary $\check D$ of the
symmetric Siegel domain
$$D:=\big\{(w,z)\in\CC^{m\times m}\times\CC^{m\times n}:(w+w^{*}-z
z^{*})\;>0\;\big\}\,,$$ which is also the interior of the convex hull
of $Q$. The inverse Cayley transform $\kappa^{-1}$, see \Ruf{UQ}, maps
the Siegel domain $D$ biholomorphically onto the bounded symmetric
domain
$$B:=\big\{(w,z)\in\CC^{m\times m}\times\CC^{m\times n}:(\One-ww^{*}-z
z^{*})\;>\;0\big\}$$ and the quadric $Q=\check D$ to an open dense part
of the (compact) \v Silov boundary 
$$\check B=\big\{(w,z)\in\CC^{m\times m}\times\CC^{m\times
n}:(ww^{*}+z z^{*})\;=\;\One\big\}$$ (called a generalized sphere
above) of $B$.  The extended quadric $\hat Q$ is CR-isomorphic to
$\check B$.

Now consider the symmetry $\sigma\in\Aut(D)$ of $D$ at the point
$e:=(\One,0)\in D$. Then $\sigma=\kappa\circ(-\id)\circ\kappa^{-1}$
with $-\id$ being the symmetry of $B$ at the origin. A simple
computation shows $\sigma(w,z)=(w^{-1},-w^{-1}z)$. But $\sigma$
extends to a symmetry of $Q$ in the sense of Section \ruf{symmetry}
and is essentially the same (up to the sign in the second variable)
as the symmetry $\gamma$ we used before. Besides these bounded
symmetric domains of type {\bf I} and their variations from Example
\ruf{EY} there are two more non-tube types of irreducible bounded
symmetric domains, all leading to symmetric standard quadrics, namely
those of types {\bf II} and {\bf V} (for this and further details in
the following see \Lit{LOSO}).

\kap{Type II} Fix an even integer $m\ge4$ and let
$j:=\kmatrix{\;\,0}{\;\;\one}{-\one}{\,0}\in\GL_{m}(\2Z)$.  Also put
$$W:=\{w\in\CC^{m\times m}:w^{J}=w\}\Steil{and}Z:=\CC^{m\times 1}\,,$$
where $w^{J}:=jw'j^{-1}$ and $w'$ is the transpose of $w$. Then $W$ is
a complex unital Jordan subalgebra of $\CC^{m\times m}$ of dimension
${m\choose2} $with $W=W^{*}$.$$Q:=\big\{(w,z)\in W\oplus
Z:w+w^{*}=zz^{*}+(zz^{*})^{J}\big\}\,,$$ is a standard quadric with
symmetry $\sigma(w,z)=(w^{-1},-w^{-1}z)$, and
$\hol(Q)\cong\so^{*}(2m{+}2)$.

\kap{Type V} Let $\OO$ be the real Cayley division algebra and
$x\mapsto \overline x\,$ its canonical (real) involution. $\OO$ is an
alternative algebra of dimension 8 over $\RR$ with $\overline x=x$ if
and only if $x\in\RR\cd1$. Denote by $\OO^{\CC}$ the complexification
of $\OO$ and extend the involution of $\OO$ to a conjugate linear
involution $z\mapsto \overline z$ of the complex alternative algebra
$\OO^{\CC}$. Then $\overline{wz}=\overline z\;\overline
w\,$ for all $w,z\in\OO^{\CC}$ and
$$Q:=\{(w,z)\in\OO^{\CC}\times\OO^{\CC}:w+\overline w=z\overline z\}$$
is a standard quadric of CR-codimension 8 in $\CC^{16}$. Again,
$\sigma(w,z)=(w^{-1},-w^{-1}z)$ is a symmetry of $Q$. Notice that for
every $w,z\in\OO^{\CC}$ with $w$ invertible there is a unital {\sl
associative} complex subalgebra of $\OO^{\CC}$ containing $w,w^{-1}$
and $z$. Furthermore, $\hol(Q)$ is isomorphic to the exceptional real
Lie algebra $\7e_{6(-14)}\,$.

\vskip7mm {\gross\noindent References} \medskip
{\klein
\parindent 15pt\advance\parskip-1pt

  }

\Ref{BAJT}Baouendi, M.S., Jacobowitz, H., Treves F.: On the analyticity of CR mappings. Ann. of Math. {\bf 122}, 365-400 (1985).
\Ref{BELO}Beloshapka, V.: On holomorphic transformations of a quadric. Math. USSR Sb. {\bf 72}, 189-205 (1992).
\Ref{ESSC}Ezhov, V., Schmalz, G.: A Matrix Poincar\'e Formula for Holomorphic Automorphisms of Quadrics of Higher Codimension. Real Associative Quadrics. J. Geom. Analysis {\bf 8}, 27-41 (1998).
\Ref{ESCH}Ezhov, V., Schmalz, G.: Automorphisms of Nondegenerate CR-Quadrics and Siegel Domains. Explicit Description. J. Geom. Analysis {\bf 11}, 441-467 (2001).
\Ref{FKAP}Fels, G., Kaup, W.: Local tube realizations of CR-manifolds and maximal abelian subalgebras. arXiv 0810.2019.
\Ref{ISKA}Isaev, A., Kaup, W.: Regularization of Local CR-Automorphisms of Real-Analytic CR-manifolds. arXiv:0906.3079.
\Ref{KAZA}Kaup W., Zaitsev, D.: On Symmetric Cauchy-Riemann Manifolds. Adv. Math. {\bf149}, 145-181 (2000).
\Ref{LANG}Lang, S: {\sl Algebra.} Graduate Texts in Mathematics 211, Springer 2005.
\Ref{LOSO}Loos, O.: {\sl Bounded symmetric domains and Jordan pairs,} Mathematical Lectures. Irvine: University of California at Irvine 1977.
\Ref{TUMA}Tumanov, A., Finite-dimensionality of the group of CR automorphisms of a standard CR manifold, and proper holomorphic mappings of Siegel domains.  Math. USSR. Izv. {\bf 32}, 655--662 (1989).

\vglue1cm\parindent0pt
Mathematisches Institut der Universit\"at\nline
72076 T\"ubingen, Auf der Morgenstelle 10\nline
GERMANY\nline
E-mail: kaup@uni-tuebingen.de

\closeout\aux\bye